\documentclass[12pt]{amsart}

\usepackage{amsmath}
\usepackage{amsfonts}
\usepackage{amssymb}

 \newtheorem{theorem}{Theorem}

\author{Sergey V. Galaev}

\title{Interior Geometry of Almost Contact K\"ahlerian Manifolds}

\begin{document}

\begin{abstract}
 In this paper, the notion  of an almost contact K\"ahlerian structure is introduced.
The interior geometry of  almost contact K\"ahlerian spaces is investigated.
On the   zero-curvature distribution of an almost contact metric structure,
as on the total space of a vector bundle, an  almost contact K\"ahlerian structure is obtained.

\end{abstract}

\maketitle

\section{Introduction}
Almost contact metric structures $(\varphi,\vec\xi,\eta,g)$
are odd-dimensional analogs of almost Hermitian structures.
There are a lot of
important interplays between these structures. The most of the
works devoted to the investigation of the geometry of manifolds
with almost contact metric structures, explicitly or not
explicitly, either use these interplays or find their
specifications. On the other hand, the presence of a smooth distribution $D$ in the geometry of almost
contact metric space gives the possibility to use the methods of the geometry of non-holonomic manifolds
   in the investigations of almost contact metric structures. Probably the
   possibility of the effective use of such approach to the
   investigation of almost contact metric spaces was stated for
   the first time in \cite{l1}. In the same time, the works, where per se an attempt to the attainment of the
   compromise on the way of the rapprochement of "holonomic" and "non-holonomic"
   points of view is done,  appeared. An example of such works is
   \cite{l2}. The main result of \cite{l2} is the construction of a new linear connection $\nabla$ on
   a contact metric space by using the Levi-Civita connection. The
   author of \cite{l2} called this connection a $D$-connection.
   This connection, in particular, satisfies the following
   property: a contact metric space is Sasakian if and only if
   $\nabla \varphi=0$ \cite[p. 1963]{l2}. The author of \cite{l2}
   writes: "As a conclusion we may say that the study of the contact distribution
   $(D, \varphi, g)$ by using the $D$-connection $\nabla$ is an
alternative to the study of the contact metric manifold $M$ via
the Levi-Civita connection" \cite[p. 1967]{l2}. The appearance of
the contact distribution $(D,\varphi,g)$ indicates the attempt to
use the methods of the non-holonomic geometry for the
investigation of almost contact metric structures. In the present
paper we use the interior connection introduced by Wagner
\cite{l3} in order to investigate the almost contact metric
structures. We develop the notion of the interior connection (we call
it a connection over a distribution), we introduce the notion of
the extended connection and generalize these notions to the
connections of the Finslerian type \cite{l4}. We show that if
$\nabla$ is an interior metric connection and $\nabla^1$ is the
corresponding extended connection, then the following statement
holds: an almost contact Hermitian space is an almost contact
K\"ahlerian space if and only if $\nabla^1\varphi=0$. The last
statement is a theorem of the proper non-holonomic geometry.

Let $(\varphi,\vec\xi,\eta,g)$ be an almost contact metric
structure (the main theses of the theory of almost contact metric
structures can be found in the excellent books \cite{l5,l6}). By
definition, an almost contact metric structure is Sasakian if it
is normal, i.e. $$N_\varphi+2d\eta\otimes\vec\xi=0,$$ where
$N_{\varphi}$ is the Nijenhuis torsion defined for the tensor
$\varphi$ and it holds $\Omega=d\eta$, where $\Omega(\vec X,\vec
Y)=g(\vec X,\varphi\vec Y)$ is the fundamental form of the
structure. Thus with an almost contact metric space we associate
 two 2-forms, $\omega=d\eta$ and $\Omega$. If these forms are equal,
we get a contact metric space, which are more simple as the more
general contact metric spaces. We will get a space with an almost
contact Hermitian structure if we refuse the condition
$\Omega=d\eta$, and the condition
$N_\varphi+2d\eta\otimes\vec\xi=0$ change to the weaker one
$$N_\varphi+2(d\eta\circ\varphi)\otimes\vec\xi=0.$$ We also do not
assume that the equality $\omega(\varphi\vec X,\varphi\vec
Y)=\omega(\vec X,\vec Y)$. Almost contact Hermitian spaces
preserve many important properties of Sasakian spaces and they
remain to be analogs of Hermitian spaces. If we make some natural
assumption about an almost contact Hermitian space, then we get an
almost contact K\'ahlerian space, these spaces are analogs of
K\"ahlerian spaces.

Following the ideology developed in the works of Schouten and
Wagner, we define the intrinsic geometry of an almost contact
metric space $X$ as the aggregate of the properties that possess
the following objects: a smooth distribution $D$ defined by a
contact form   $\eta$; an admissible field of endomorphisms
$\varphi$ of $D$ (which we call an admissible almost complex
structure) satisfying $\varphi^{2}=-1$; an admissible Riemannian
metric field $g$ that is related to $\varphi$ by
$g(\varphi\vec{X},\varphi\vec{Y})=g(\vec{X},\vec{Y})$, where
$\vec{X}$ and $\vec{Y}$ are admissible vector fields. To the
objects of the intrinsic geometry of an almost contact metric
space one should ascribe also the  objects derived from the just
mentioned: the 2-form $\omega=d\eta$; the vector field $\vec{\xi}$
(which is called the Reeb vector field) defining the closing
$D^{\bot}$ of $D$, i.e. $\vec{\xi} \in D^{\bot}$, and given by the
equalities $\eta(\vec{\xi})=1$, ${\rm ker\,} \omega={\rm span}
(\vec{\xi})$ in the case when the 2-form  $\omega$ is of maximal
rank; the intrinsic connection $\nabla$ that  defines the parallel
transport of admissible vectors along admissible curves and is
defined by the metric $g$; the connection $\nabla^{1}$ that is a
natural extension of the connection $\nabla$ which accomplishes
the parallel transport of admissible vectors along arbitrary
curves of the manifold $X$.

Besides the introduction, the paper contains  4 sections. In
Section \ref{sec2} we introduce the notion of an almost contact
K\"ahlerian structure. In Section \ref{sec3} we discuss
connections over a distribution and the extended connections. The
notion of the connection over a distribution were known before
(see e.g. \cite{l7}). The connection over a distribution was used
in the geometry of contact structures in \cite{l4,l8}. In these
works also the notion of the extended connection was used. Per se,
the extended connection was defined for the first time by Wagner
in a little bit another context and in another terms in \cite{l3}
in order to construct the curvature tensor of a non-holonomic
manifold. In Section \ref{sec4} we give the main results about the
interior geometry of almost contact K\"ahlerian spaces. This
section contains the main results of the paper: an almost contact
Hermitian structure is an almost contact K\"ahlerian structure if
and only if $\nabla^1\varphi=0$, where $\nabla^1$  is the interior
metric torsion-free connection. The last sections contains an
example of an almost contact K\"ahlerian space that is not a
Sasakian space.

\section{Almost contact K\"ahlerian structure}\label{sec2}

Let $X$  be a smooth manifold of  an odd dimension $n$, $n\geq 3$.
Denote by $\Xi(X)$ the $C^{\infty}(X)$-module of smooth vector
fields on $X$. All manifolds, tensors and other geometric objects
will be assumed to be smooth of the class $C^{\infty}$.  An almost
contact metric structure on
 $X$ is an aggregate
 $(\varphi, \vec{\xi}, \eta, g)$ of  tensor fields on $X$, where $\varphi$ is a tensor field of type $(1, 1)$, which is called the structure
 endomorphism, $\vec{\xi}$ and $\eta$  are a vector and a covector,
 which are called the structure vector and the contact form,
 respectively, and  $g$ is a (pseudo-)Riemannian metric.
 Moreover,
$$\eta(\vec{\xi})=1,\quad \varphi(\vec{\xi})=0,\quad \eta \circ
\varphi=0,$$
$$\varphi^2\vec{X}=-\vec{X}+\eta(\vec{X})\vec{\xi},\quad
g(\varphi\vec{X},\varphi\vec{Y})=g(\vec{X},\vec{Y})-\eta(\vec{X})\eta(\vec{Y})$$
for all $\vec{X}, \vec{Y} \in \Xi(X)$. The skew-symmetric tensor
 $\Omega(\vec{X}, \vec{Y})=g(\vec{X}
\varphi\vec{Y})$ is  called the fundamental tensor of the
structure. A manifold with a fixed almost contact metric structure
is called an almost contact metric manifold. If
 $\Omega=d\eta$ holds, then the almost contact metric structure is called a contact metric structure.
An almost contact metric structure is called normal if
$$N_{\varphi}+2d\eta\otimes\vec{\xi}=0,$$ where $N_{\varphi}$ is the
Nijenhuis torsion defined for the tensor $\varphi$. A normal
contact metric structure is called a Sasakian structure. A
manifold with a given Sasakian structure is called a Sasakian
manifold. Let $D$ be the smooth distribution of codimension 1
defined by the form  $\eta$, and $D^\bot={\rm span}(\vec{\xi})$ be
the closing of $D$. If  the restriction of the 2-form
$\omega=d\eta$ to the distribution $D$  is non-degenerate, then
the vector $\vec{\xi}$ is uniquely defined by the condition
$$\eta(\vec{\xi})=1,\quad {\rm ker\,} \omega={\rm span}
(\vec{\xi}),$$ and  it is called the Reeb vector field.

We say that an almost contact metric structure is almost normal,
if it holds
\begin{equation}\label{q1}N_{\varphi}+2d\eta\otimes\vec{\xi}=0.\end{equation}
In what follows, an almost normal almost contact metric space will
be called {\it an almost contact Hermitian space}. An almost
contact Hermitian space is called {\it an almost contact
K\"ahlerian space}, if its fundamental form is closed. The
following obvious theorem shows  the difference between  a
normal almost contact metric structure and an almost contact
Hermitian structure.

\begin{theorem}\label{thN1} An almost contact Hermitian structure is normal if
and only if it holds $$\omega(\varphi\vec u,\varphi\vec
v)=\omega(\vec u,\vec v),\quad\vec u,\vec v\in\Gamma
D.$$\end{theorem}

It is obvious that an almost normal contact metric structure is a
Sasakian structure. Sasakian manifolds are popular among the
researchers of almost contact metric spaces by the following two
reasons. On one hand, there exist a big number of interesting and
deep examples of Sasakian structures (see e.g. \cite{l6}), on the
other hand, the Sasakian manifolds have very important and natural
properties.

We say that a coordinate map $K(x^\alpha)$
$(\alpha,\beta,\gamma=1,...,n)$ $(a,b,c,e=1,...,n-1)$ on a
manifold $X$ is adapted to the non-holonomic manifold  $D$ if
$$D^{\bot}={\rm span}\left(\frac{\partial}{\partial
x^{n}}\right)$$ holds \cite{l1}.

Let $P:TX\rightarrow D$ be the projection map defined by the
decomposition $TX=D\oplus D^{\bot}$ and let $K(x^{\alpha})$ be an
adapted coordinate map.  Vector fields
$$P(\partial_{a})=\vec{e}_{a}=\partial_{a}-\Gamma^{n}_{a}\partial_{n}$$
are linearly independent, and linearly generate the system $D$
over the domain of the definition of the coordinate map: $$D={\rm
span}(\vec{e}_{a}).$$ Thus we have on $X$ the non-holonomic field
of bases  $(\vec{e}_{a},\partial_{n})$ and the corresponding field
of cobases $$(dx^a,\theta^{n}=dx^{n}+\Gamma^{n}_{a}dx^{a}).$$
 It can be checked directly that
$$[\vec{e}_{a},\vec{e}_{b}]=M^{n}_{ab}\partial_{n},$$ where the
components  $M^{n}_{ab}$ form the so-called tensor of
non-holonomicity \cite{l3}. Under assumption that for all adapted
coordinate systems it holds  $\vec{\xi}=\partial_{n}$, the
following equality takes place
$$[\vec{e}_{a},\vec{e}_{b}]=2\omega_{ba}\partial_{n},$$ where
$\omega=d\eta$.  We say also that the basis
$$\vec{e}_{a}=\partial_{a}-\Gamma^{n}_{a}\partial_{n}$$ is adapted,
as the basis defined by an adapted coordinate map. Note that
$\partial_n\Gamma^n_a=0$.

We call a tensor field defined on an almost contact metric
manifold admissible (to the distribution $D$) if it vanishes
whenever its vectorial argument belongs to the closing
 $D^\bot$ and its covectorial argument is proportional to the form  $\eta$.
 The coordinate form of an admissible tensor field
of type  $(p,q)$ in an adapted coordinate map looks like
$$
t=t^{a_{1},...,a_{p}}_{b_{1},...,b_{q}}\vec{e}_{a_{1}}\otimes...\otimes\vec{e}_{a_{p}}\otimes
dx^{b_{1}}\otimes...\otimes dx^{b_{q}}.
$$
In particular, an admissible vector field is a vector field
tangent to the distribution $D$, and an admissible 1-form is a
1-form that is zero on the closing  $D^\bot$. It is clear that any tensor
structure defined on the manifold  $X$ defines on it a unique
admissible tensor structure of the same type. From the definition
of an almost contact structure it follows that the field of
endomorphisms $\varphi$ is an admissible tensor field of type $(1,
1)$. The field of endomorphisms  $\varphi$ we call an admissible
almost complex structure, taking into the account its properties.
The 2-form $\omega=d\eta$ is also an admissible tensor field and
it is natural to call it an admissible symplectic form.

All constructions done by Wagner in \cite{l3} are grounded on the
usage of adapted coordinates (Wagner called such coordinates by
gradient coordinates). Adapted coordinates are used in the
foliation theory \cite{l9}. It seems that in the theory of almost
contact metric spaces, the adapted coordinates were used in essence
only in the works \cite{l1,l2,l4}.

One of the main notions of this work is the notion of an
admissible integrable tensor structure. In the definition of an
admissible integrable tensor structure, the words "integrable"
and "admissible" should be consider in the semantic union. We call
an admissible tensor field {\it integrable} if there is an open
neighborhood of each point of the manifold $X$ and admissible
coordinates on it such that the components  of the tensor fields
are constant with respect to these coordinates. The form
$\omega=d\eta$ is an example of an admissible tensor structure. If
the distribution $D$ is integrable, then any admissible integrable
structure is an integrable structure on the manifold $X$. The
following facts show that the notion of an integrable admissible
tensor structure is natural. As it is known, the integrable
closing $D^\bot$ defines a foliation with one-dimensional lives.
If one defines on this foliation a structure of a smooth
manifold, then that any integrable tensor structure defines on
this manifold an integrable tensor structure in the usual sense.
Below we consider some important ideas of the development of the
notion of the integrability in the geometry of almost contact
metric structures. The following two theorems show the meaning of
the notion of an integrable tensor structure in the context of our
investigations.

\begin{theorem}\label{Th2}  The admissible almost complex strructure
  $\varphi$ is integrable if and only if $P(N_{\varphi})=0$ holds.
\end{theorem}

\textbf{Proof.} The expression of the non-zero components of the
Nijenhuis torsion tensor
$$N_{\varphi}(\vec{X}, \vec{Y})=[\varphi \vec{X}, \varphi
\vec{Y}]+\varphi^{2}[\vec{X},\vec{Y}]-\varphi[\varphi\vec{X},\vec{Y}]-\varphi[\vec{X},\varphi\vec{Y}]$$
of the tensor $\varphi$ in adapted coordinates has the form:
\begin{align}\label{eq2.a} N^{e}_{ab}&=\varphi^{c}_{a}\vec{e}_{c}\varphi^{e}_{b}-
\varphi^{c}_{b}\vec{e}_{c}\varphi^{e}_{a}+\varphi^{e}_{c}\vec{e}_{b}\varphi^{c}_{a}-\varphi^{e}_{d}\vec{e}_{a}\varphi^{d}_{b},\\
\label{eq2.b}
N^{n}_{ab}&=2\varphi^{c}_{a}\varphi^{d}_{b}\omega_{dc}, \\
\label{eq2.c}
N^{e}_{na}&=-\varphi^{e}_{c}\partial_{n}\varphi^{c}_{a}.\end{align}
Thus the equality $P(N_{\varphi})=0$ is equivalent to the
condition that \eqref{eq2.a} and \eqref{eq2.c} are zero.

Conversaly, suppose that $P(N_{\varphi})=0$. Consider a
sufficiently small neighborhood  $U$  of an arbitrary point of the
manifold $X$. Assume that  $U=U_{1} \times U_{2}$, $TU={\rm
span}(\partial_{a}) \oplus {\rm span}(\partial_{n})$. We set the
natural denotation $T(U_{1})={\rm span}(\partial_{a})$. We define
over the set   $U$ the isomorphism of bundles $\psi:D\rightarrow
T(U_{1})$ by the formula $\psi(\vec{e}_{a})=\partial_{a}$. This
endomorphism induces an almost complex structure on the manifold
$U_{1}$. This complex structure is integrable due to the equality
$P(N_{\varphi})=0$. Indeed, from \eqref{eq2.c} it follows that the
right hand side part of \eqref{eq2.a} coincides with the torsion
of the almost complex structure induced on the manifold $U_{1}$.
Choosing an appropriate coordinate system on $U_{1}$, and
consequently, an appropriate adapted coordinate system on the
manifold  $X$, we get a coordinate map with respect to that the
components of the endomorphism  field $\varphi$ are constant.
$\Box$
\begin{theorem}\label{Th3} An almost contact metric structure is
is an almost contact Hermitian structure if and only if the
admissible almost complex structure $\varphi$ is integrable.
\end{theorem}

{\bf Proof.} Equality \eqref{q1} written in adapted coordinates is
equivalent to the condition that the right hand sides of
\eqref{eq2.a} and \eqref{eq2.c} are zero. This and Theorem
\ref{Th2} gives the proof of the theorem. $\Box$

Note that we have de facto proved the equality
$$P(N_\varphi)=N_\varphi+2(d\eta\circ\varphi)\otimes\vec\xi.$$

Using adapted coordinates we introduce the following admissible
tensor fields: $$h^a_b=\frac{1}{2}\partial_n\varphi^a_b,\quad
C_{ab}=\frac{1}{2}\partial_ng_{ab},\quad C^a_b=g^{da}C_{db},\quad
\psi^b_a=g^{db}\omega_{da}.$$ We denote by $\tilde\nabla$ and
$\tilde\Gamma^\alpha_{\beta\gamma}$ the Levi-Civita connection and
the Christoffel symbols of the metric $g$. The proof of the
following theorem follows from direct computations.

\begin{theorem}\label{Th4} The Christoffel symbols of the Levi-Civita connection
of an almost contact metric space with respect to adapted
coordinates are the following:
$$\tilde\Gamma^c_{ab}=\Gamma^c_{ab},\quad \tilde
\Gamma^n_{ab}=\omega_{ba}-C_{ab},\quad \tilde \Gamma^b_{an}=\tilde
\Gamma^b_{na}=C^b_a-\psi^b_a,\quad \tilde \Gamma^n_{na}=\tilde
\Gamma^a_{nn}=0,$$ where $$\Gamma^a_{bc}=\frac{1}{2}g^{ad}(\vec e_b g_{cd}-\vec e_c g_{bd}-\vec e_d g_{bc}).$$
\end{theorem}

In the case of a contact metric space the Christoffel symbols
of the Levi-Civita connection are found in \cite{l2}.

\section{Connection over a distribution. The extended connection}\label{sec3}

An intrinsic linear connection on a non-holonomic manifold $D$ is
defined in \cite{l3} as a map  $$\nabla:\Gamma D \times \Gamma
D \rightarrow \Gamma D $$ that satisfies the following conditions:
\begin{align*}
1)\quad &
\nabla_{f_1\vec{u}_1+f_2\vec{u}_2}=f_1\nabla_{\vec{u}_1}+f_2\nabla_{\vec{u}_2};\\
2)\quad &
\nabla_{\vec{u}}f\vec{v}=f\nabla_{\vec{u}}\vec{v}+(\vec{u}f)\vec{v},
\end{align*} where $\Gamma D$ is the module of admissible vector fields. The
Christoffel symbols are defined by the relation $$
\nabla_{\vec{e}_{a}}\vec{e}_{b}=\Gamma^{c}_{ab}\vec{e}_{c}. $$

The torsion $S$ of the intrinsic linear connection is defined by
the formula
$$
S(\vec{X},\vec{Y})=\nabla_{\vec{X}}\vec{Y}-\nabla_{\vec{Y}}\vec{X}-p[\vec{X},\vec{Y}].
$$ Thus with respect to an adapted coordinate system it holds  $$
S^{c}_{ab}=\Gamma^{c}_{ab}-\Gamma^{c}_{ba}. $$ 

The action of an interior linear connection can be extended in a natural way to
admissible tensor fields. An important example of an interior linear connection is
the interior metric connection that is uniquely defined by the conditions $\nabla g=0$ and 
$S=0$ \cite{l4}. With respect to the adapted coordinates it holds 
\begin{equation}\label{for3} \Gamma^a_{bc}=\frac{1}{2}g^{ad}(\vec e_b g_{cd}-\vec e_c g_{bd}-\vec e_d g_{bc}).
\end{equation}

 In the same way as
a linear connection on a smooth manifold, an intrinsic connection
can be defined by giving a horizontal distribution over the total
space of some vector bundle. The role of such bundle plays the
distribution $D$. The notion of   {\it a connection over a
distribution} was applied 
to non-holonomic manifolds with admissible Finsler metrics in
\cite{l4,l8}. One says that over a distribution  $D$ a
connection is given if the distribution
$\tilde{D}=\pi^{-1}_{*}(D)$, where $\pi:D \rightarrow X$ is the
natural projection, can be decomposed into a direct some of the
form $$\tilde{D}=HD \oplus VD,$$ where $VD$ is the vertical
distribution on the total space $D$.

Let us introduce a structure of a smooth manifold on $D$. This
structure is defined in the following way. To each adapted
coordinate map  $K(x^\alpha)$ on the manifold  $X$ we put in
correspondence the coordinate map
$\tilde{K}(x^{\alpha},x^{n+\alpha})$ on the manifold  $D$, where
$x^{n+\alpha}$ are the coordinates of an admissible vector with
respect to the basis
$$\vec{e}_{a}=\partial_{a}-\Gamma^{n}_{a}\partial_{n}.$$ The
constructed over-coordinate map will be called adapted. Thus the
assignment of a connection over a distribution is equivalent to
the assignment of the object $G^{a}_{b}(X^{a},X^{n+a})$  such that
$$HD={\rm span}(\vec{\epsilon}_{a}),$$ where
$\vec{\epsilon}_{a}=\partial_{a}-\Gamma^{n}_{a}\partial_{n}-G^{b}_{a}\partial_{n+b}$.
If it holds
$$G^{a}_{b}(x^{a},x^{n+a})=\Gamma^{a}_{bc}(x^{a})x^{n+c},$$
then  the connection over the distribution $D$  is defined by  the
interior linear connection. In \cite{l4} the notion of the
prolonged connection was introduced. The prolonged
 connection can be obtained from an intrinsic connection by the equality
$$TD=\tilde{HD} \oplus VD,$$ where $HD \subset \tilde{HD}$.
 Essentially, the prolonged connection is a connection in a vector bundle.
As it follows from the definition of the extended connection, for
its assignment (under the condition that a connection on the
distribution is already defined) it is enough to define a vector
field on the manifold $D$ that has the following coordinate form:
$\vec u=\partial_n-G^a_n\partial_{n+a}$. The components of the
object $G^a_n$ are transformed as the components of a vector on
the base. Setting $G^a_n=0$, we get an extended connection denoted by $\nabla^1$.
In \cite{l3} the admissible tensor field
$$R(\vec u,\vec v)\vec w=\nabla_{\vec u}\nabla_{\vec v}\vec w-\nabla_{\vec v}\nabla_{\vec u}\vec w-
\nabla_{p[\vec u,\vec v]}\vec w$$ is called by Wagner the first Schouten curvature tensor. With respect to the adapted coordinates it holds $$R^a_{bcd}=2\vec e_{[a}\Gamma^d_{b]c}+2\Gamma^d_{[a||e||}\Gamma^e_{b]c}.$$
If the distribution $D$ does not contain any integrable subdistribution of dimension $n-2$, then the Schouten 
curvature tensor  is zero if and only if the parallel transport of admissible vectors does
not depend on the curve \cite{l3}. We say that the Schouten tensor is the curvature tensor of the
distribution $D$. If this tensor is zero, we say that the distribution $D$ is a zero-curvature distribution.
Note that the partial derivatives are components of an admissible tensor field \cite{l3}.

\section{Properties of almost contact K\"ahlerian spaces related to the usage of the cobnnection over a distribution}\label{sec4}

Let $(\varphi,\vec\xi ,\eta ,g)$ be an almost contact metric structure. In \cite{l1}, the following theorem is proved:

\begin{theorem}\label{Th5} Let $\nabla$ be a torsion-free interior linear connection on an almost contact metric space $X$. Then there exists on $X$ a connection with the torsion
$$S(\vec{x},\vec{y})=\frac{1}{4}P\left(N_{\varphi }\right)\left(\vec{x},\vec{y}\right),\quad \vec{x},\vec{y}\in \Gamma D$$ and compatible with $\varphi $. \end{theorem}

The following theorem is the corollary of Theorem \ref{Th5}.

\begin{theorem}\label{Th6} An almost contact metric space admits a torsion-free interior connection $\nabla$
such that ${\nabla }^1\varphi =0$ if and only if the admissible structure $\varphi$ is integrable.
\end{theorem}

{\bf Proof.} Let 
$\nabla$ be a torsion-free connection such that ${\nabla }^1\varphi =0$. Applying this
 $\nabla$ to the proof of Theorem \ref{Th5}, we get $$S(\vec{x},\vec{y})=\frac{1}{4}P\left(N_{\varphi }\right)\left(\vec{x},\vec{y}\right)=0,\quad \vec{x},\vec{y}\in \Gamma D.$$ Adding to this condition the equality ${\partial }_n{\varphi }^a_b=0$, we get $$P\left(N_{\varphi }\right)\left(\vec{x},\vec{y}\right)=0,\quad \vec{x},\vec{y}\in TX.$$ By Theorem \ref{Th2}, this is equivalent to the integrability of $\varphi $. The converse statement is obvious. $\Box$

\begin{theorem}\label{Th7} An almost contact metric structure is an almost contact K\"ahlerian structure if and only if ${\nabla }^1\varphi =0$, where $\nabla $ is the interior torsion-free metric connection.
\end{theorem}

{\bf Proof.} According \cite{l2}, any almost contact metric space satisfies  the following equality:
\begin{multline}\label{q4}2g(({\widetilde{\nabla }}_{\vec{x}}\varphi )\vec{y},\vec{z})=3d\Omega \left(\vec{x},\varphi \vec{y},\varphi \vec{z}\right)-3d\Omega \left(\vec{x},\vec{y},\vec{z}\right)+g\left(N^{\left(1\right)}\left(\vec{y},\vec{z}\right),\varphi \vec{x}\right)\\+N^{\left(2\right)}\left(\vec{y},\vec{z}\right)\eta \left(\vec{x}\right)+
2d\eta \left(\varphi \vec{y},\vec{x}\right)\eta \left(\vec{z}\right)-2d\eta \left(\varphi \vec{z},\vec{x}\right)\eta (\vec{y}),\end{multline} where
 $$N^{(1)}=N_{\varphi }+2d\eta \otimes \vec{\xi },\quad N^{\left(2\right)}(\vec{x},\vec{y})=\left(L_{\varphi \vec{x}}\eta \right)\vec{y}-\left(L_{\varphi y}\eta \right)\vec{x}.$$
Theorem \ref{Th2} and the definition of an almost contact K\"ahlerian structure allow us to assume in what follows that the almost contact metric structure $(\varphi,\vec\xi,\eta, g)$ is almost normal. In this case, $$P(N_{\varphi })=N_{\varphi }+2(d\eta \circ \varphi )\otimes \vec{\xi }=0.$$  Thus,  
$$N^{(1)}=2(d\eta \otimes \vec{\xi }-(d\eta \circ \varphi )\otimes \vec{\xi }),$$ and the equality
\eqref{q4} takes the simpler form:
\begin{multline}\label{q5} 2g(({\widetilde{\nabla }}_{\vec{x}}\varphi )\vec{y},\vec{z})=3d\Omega \left(\vec{x},\varphi \vec{y},\varphi \vec{z}\right)-3d\Omega \left(\vec{x},\vec{y},\vec{z}\right)+N^{\left(2\right)}\left(\vec{y},\vec{z}\right)\eta \left(\vec{x}\right)\\+2d\eta \left(\varphi \vec{y},\vec{x}\right)\eta \left(\vec{z}\right)-2d\eta \left(\varphi \vec{z},\vec{x}\right)\eta (\vec{y}).\end{multline} 
{\it Sufficiency.} Substituting to \eqref{q5} first  $\vec{x}={\vec{e}}_a$, $\vec{y}={\partial }_n$, $\vec{z}={\vec{e}}_c$, and then $\vec{x}={\vec{e}}_a$, $\vec{y}={\vec{e}}_b$, $\vec{z}={\vec{e}}_c$, we get $d{\Omega }_{abn}=0$ and $d{\Omega }_{abc}=0$, respectively. This means that $d\Omega =0$.

{\it Necessity.} Suppose that $d\Omega =0$. We may rewrite \eqref{q5} in the form 
\begin{equation}\label{q6}
2g(({\widetilde{\nabla }}_{\vec{x}}\varphi )\vec{y},\vec{z})=N^{\left(2\right)}\left(\vec{y},\vec{z}\right)\eta \left(\vec{x}\right)+2d\eta \left(\varphi \vec{y},\vec{x}\right)\eta \left(\vec{z}\right)-2d\eta \left(\varphi \vec{z},\vec{x}\right)\eta (\vec{y}).  \end{equation}
Substituting  $\vec{x}={\vec{e}}_a$, $\vec{y}={\vec{e}}_b$, $\vec{z}={\vec{e}}_c$ to \eqref{q6}, we get ${\nabla }_a{\varphi }^b_c=0$. $\Box$

In the rest of this section we formulate and prove a theorem generalizing the following classical result \cite{l2}: an almost contact metric space is a Sasakian space if and only if the following equality holds:
                           $$({\widetilde{\nabla }}_{\vec{x}}\varphi )\vec{y}=g(\vec{x},\vec{y})\vec{\xi }-\eta \left(\vec{y}\right)\vec{x}.$$                                        

\begin{theorem}\label{Th8} An almost contact Hermitian structure is an almost contact K\"ahlerian structure if and only if it holds
\begin{equation} \label{q7} 
({\widetilde{\nabla }}_{\vec{x}}\varphi )\vec{y}=d\eta \left(\varphi \vec{y},\vec{x}\right)\vec{\xi }+\eta \left(\vec{y}\right)\left(\varphi \circ \psi \right)\left(\vec{x}\right)-\eta \left(\vec{x}\right)\left(\varphi \circ \psi -\psi \circ \varphi \right)\vec{y}. 
\end{equation} \end{theorem}
{\bf Proof.} The equality \eqref{q7} is equivalent to the following conditions:
 $$\nabla \varphi =0,\quad {\partial }_n{\varphi }^a_b=0,\quad {\partial }_ng_{ab}=0.$$
 The last two equalities are written with respect to the adapted coordinates. The first two equalities can be unified by the condition ${\nabla }^1\varphi =0$, which implies  ${\partial }_ng_{ab}=0$. $\Box$

\section{Almost contact metric structures over a zero-curvature distribution}\label{sec5}

Consider the vector bundle $(D, \pi , X)$, where $D$ is the distribution of the contact metric structure
$(\varphi,\vec\xi ,\eta , g)$. If the distribution $D$ is a zero-curvature distribution and it does not
contain any involutive subdistribution of dimension $n-2$, then the equality $P^a_{bc}=0$ holds \cite{l3}.
In what follows we assume that $n>3$. On the total space $D$ of the vector bundle under the consideration we define an almost contact metric structure $(\tilde D,\tilde g,J,d(\pi^*\circ \eta),D)$ by setting
 $$\tilde g(\vec\epsilon_a,\vec\epsilon_b)=\tilde g(\partial_{n+a},\partial_{n+b}=\tilde g(\vec e_a,\vec e_b,\quad \tilde g(\vec\epsilon_a,\partial_n)=\tilde g(\partial_{n+a},\partial_n=0,$$  
 $$J(\vec\epsilon_a)=\partial_{n+a},\quad J(\partial_{n+a}=-\vec\epsilon_a,\quad J(\partial_n)=0,$$ 
$\tilde D=\pi_*^{-1}(D)$, $$\tilde{D}=HD\oplus VD,$$   $VD$ is the vertical distribution on the total space $D$, and $HD$
is the horizontal space defined by the interior linear connection.
 The vector fields $$\vec\epsilon_a=\partial_a-\Gamma^n_a\partial_n-\Gamma^b_{ac}x^{n+c}\partial_{n+b},\quad\partial_n,\quad\partial_{n+a}$$ define on $D$ a non-holonomic field of bases, and the forms
 $$dx^a,\quad\Theta^n=dx^n+\Gamma^n_adx^a,\quad \Theta^{n+a}=dx^{n+a}+\Gamma^a_{bc}x^{n+b}dx^c$$ define the corresponding field of cobases. The vector field $\partial_n$ is the the Reeb vector field of the almost contact metric structure 
$(\tilde D,\tilde g,J,d(\pi^*\circ \eta))$.

\begin{theorem}\label{Th9} The almost contact metric structure $(\tilde D,\tilde g,J,d(\pi^*\circ \eta))$ is an almost contact metric structure if and only if the distribution $D$ is a zero-curvature distribution. \end{theorem}
{\bf Proof.} It is easy to check that the following holds: 
\begin{align} \label{q8}
\left[{\vec{\varepsilon }}_a,{\vec{\varepsilon }}_b\right]&=2{\omega }_{ba}{\partial }_n+R^e_{abc}x^{n+c}{\partial }_{n+e},\\
\label{q9}
[{\vec{\varepsilon }}_a,{\partial }_n]&=x^{n+c}P^b_{ac}{\partial }_{n+b},\\
\label{q10}
\left[{\vec{\varepsilon }}_a,{\partial }_{n+b}\right]&={\Gamma }^c_{ab}{\partial }_{n+c}.
\end{align}
These equalities directly imply 
\begin{align*}
N_J\left({\vec{\varepsilon }}_a,{\vec{\varepsilon }}_b\right)&=-R^e_{abc}x^{n+c}{\partial }_{n+e},\\
N_J\left({\partial }_{n+a},{\partial }_{n+b}\right)&=2{\omega }_{ba}{\partial }_n+R^e_{abc}x^{n+c}{\partial }_{n+e},\\
N_J\left({\vec{\varepsilon }}_a,{\partial }_{n+b}\right)&=0,\\
N_J\left({\vec{\varepsilon }}_a,{\partial }_n\right)&=
N_J\left({\partial }_{n+a},{\partial }_n\right)=-x^{n+c}P^b_{ac}{\partial }_{n+b}.\end{align*}
These equalities yield  the proof of the theorem. $\Box$

Let us show that the structure $(\tilde D,\tilde g,J,d(\pi^*\circ \eta))$ is not normal. 
It holds
$$N_J\left({\partial }_{n+a},{\partial }_{n+b}\right)+2d\widetilde{\eta }\left({\partial }_{n+a},{\partial }_{n+b}\right){\partial }_n=2{\omega }_{ba}{\partial }_n+R^e_{abc}x^{n+c}{\partial }_{n+e}.$$
It is clear that this expression can not be zero.

\begin{theorem}\label{Th10} The almost contact metric structure $(\tilde D,\tilde g,J,d(\pi^*\circ \eta))$ is an almost contact K\"ahlertian structure if and only if $(\varphi,\vec\xi,\eta,g)$ is a Sasakian structure with the zero-curvature distribution.\end{theorem}

{\bf Proof.} It can be checked directly that $d\Omega=0$ if and only if $d\tilde \Omega=0$, where $\tilde \Omega(\vec x,\vec y)=g(\vec x,J\vec y)$. This proves the theorem. $\Box$

Almost contact metric spaces of zero interior curvature appear in mechanics and physics.

Wagner \cite{l3,l10} paid a big attention to non-holonomic manifolds of zero curvature. In particular, in \cite{l10}, Wagner defines a non-holonomic manifold of zero curvature that is a geometric model of a solid body under non-holonomic constrains.

\bibliographystyle{unsrt}

\begin{thebibliography}{90}

\bibitem{l1} S.V. Galaev, {\it The intrinsic geometry of almost contact metric
manifolds}, Izv. Saratov. Univ. Mat. Mekh. Inform., 12:1 (2012),
16--22.

\bibitem{l2} A. Bejancu, {\it K\"hler contact distributions.} J. Geom. Phys. 60
(2010), no. 12, 1958–-1967.

\bibitem{l3} V.V.~Wagner, {\it Geometry of  $(n - 1)$-dimensional nonholonomic manifold
 in an  $n$-dimensional space},  Proc. Sem. on vect. and tens. analysis (Moscow
 Univ.).  5 (1941), 173--255.

\bibitem{l4}
A. V. Bukusheva, S. V. Galaev, {\it Almost contact metric
structures defined by connection over distribution with admissible
Finslerian metric}, Izv. Saratov. Univ. Mat. Mekh. Inform., 12:3
(2012), 17–22

\bibitem{l5} G. Pitis, {\it Geometry of Kenmotsu manifolds.} Publishing House
of Transilvania University of Brasov, Brasov, 2007. iv+160 pp.

\bibitem{l6} D.E.~Blair, {\it Contact manifolds in Riemannian geometry},  Berlin-New York: Springer-Verlag, 1976, 146~p.

\bibitem{l7} A. Vershik, L. Faddeev, {\it Differencial geometry
and Lagrangian mechanics with constrains} DAN USSR 202 (3),
555--557 (1972).

\bibitem{l8} A. V. Bukusheva, {\it About the geometry of
foliations on distributions with Finslerian metrics}, Proc. Penz.
Uni. Phys. Math. and Tech. Sci. 30 (2012), 33--38.

\bibitem{l9} M.A.~Malakhal'tsev, {\it Foliations with leaf structures}, Geometry, 7. J. Math. Sci. (New York) 108 (2002), no. 2,
188--210.

\bibitem{l10} V.V. Wagner, {\it Geometric interpretation of the motion of nonholonomic dynamical systems.} Abh. Sem. Vektor- und Tensoranalysis [Trudy Sem. Vektor. Tenzor. Analizu] 5, (1941) 301--327.


\end{thebibliography}

\vskip2cm

Saratov State University,\\
Chair of Geometry\\
E-mail: sgalaev(at)mail.ru

\end{document}